\documentclass[reqno,12pt]{amsart}

\usepackage{color} 
\usepackage{ifpdf}
\ifpdf
    \usepackage[pdftex]{graphicx}
    \usepackage[pdftex]{hyperref}
    \hypersetup{
        unicode=false,          
        pdftoolbar=true,        
        pdfmenubar=true,        
        pdffitwindow=false,     
        pdfstartview={FitH},    
        pdftitle={MCP Article},      
        pdfauthor={Michael Holst},   
        pdfsubject={Mathematics},    
        pdfcreator={Michael Holst},  
        pdfproducer={Michael Holst}, 
        pdfkeywords={PDE, analysis, mathematical physics}, 
        pdfnewwindow=true,      
        colorlinks=true,        
        linkcolor=red,          
        citecolor=blue,         
        filecolor=magenta,      
        urlcolor=cyan           
    }

\else
    \usepackage{graphicx}
    \usepackage{pstricks}
    
    \newcommand{\href}[2]{#2}
\fi

\usepackage{times}
\usepackage{amsfonts}

\usepackage{amsmath}
\usepackage{amsthm}
\usepackage{amssymb}
\usepackage{amsbsy}
\usepackage{amscd}

\usepackage{enumerate}
\usepackage{verbatim}
\usepackage{subfigure}

\newtheorem{theorem}{Theorem}[section]

\newtheorem{algorithm}[theorem]{Algorithm}

\numberwithin{equation}{section}  

  \newcounter{mnote}
  \let\oldmarginpar\marginpar
    \renewcommand\marginpar[1]{\-\oldmarginpar[\raggedleft\footnotesize #1]%
    {\raggedright\footnotesize #1}}

\definecolor{myblue}{rgb}{0.2,0.2,0.7}
\definecolor{mygreen}{rgb}{0,0.6,0}
\definecolor{mycyan}{rgb}{0,0.6,0.6}
\definecolor{myred}{rgb}{0.9,0.2,0.2}
\definecolor{mymagenta}{rgb}{0.9,0.2,0.9}
\definecolor{mywhite}{rgb}{1.0,1.0,1.0}
\definecolor{myblack}{rgb}{0.0,0.0,0.0}

\newcommand{\leqs}{\leqslant}      
     
\def\FETK{{\sc FETK}}
\newcommand{\semr}[1]{\ensuremath{\mbox{#1}}}
\def\jump#1{\bigg[\!\!\bigg[#1\bigg]\!\!\bigg]}
\def\sjump#1{\big[\!\big[#1\big]\!\big]}
\newcommand{\tpiper}{\ensuremath{\right|\!\right|\!\right|}}
\newcommand{\tpipel}{\ensuremath{\left|\!\left|\!\left|}}
\newcommand{\energy}[1]{\ensuremath{\tpipel #1 \tpiper}}
\newcommand{\calT}{\mathcal{T}}
\newcommand{\diff}{\ensuremath{\,d}}

\begin{document}

\title[Goal-Oriented Adaptivity and Multilevel Preconditioning
       for the PBE]
      {Goal-Oriented Adaptivity and Multilevel Preconditioning
       for the Poisson-Boltzmann Equation}

\author[B. Aksoylu]{Burak Aksoylu$^1$}
\thanks{$^1$
TOBB University of Economics and Technology, Department 
of Mathematics, Ankara, 06560, Turkey and
Louisiana State University, Department of Mathematics, 
Baton Rouge, LA 70803, USA.}
\email{baksoylu@etu.edu.tr}

\author[S. Bond]{Stephen~D. Bond$^2$}
\thanks{$^2$
Multiphysics Simulation Technologies Department,
Sandia National Laboratories, Albuquerque, NM 87185.}
\email{sdbond@sandia.gov}

\author[E. Cyr]{Eric~C. Cyr$^3$}
\thanks{$^3$
Scalable Algorithms Department,
Sandia National Laboratories, Albuquerque, NM 87185.}
\email{eccyr@sandia.gov}

\author[M. Holst]{Michael Holst$^4$}
\thanks{$^4$
Department of Mathematics,
University of California San Diego,
La Jolla, CA 92093.}
\email{mholst@math.ucsd.edu}

\date{\today}

\keywords{Poisson-Boltzmann equation, adaptive finite element methods, multilevel preconditioning, goal-oriented {\em a posteriori} error estimation, solvation free energy, electrostatics}

\begin{abstract}
  In this article, we develop goal-oriented error indicators to drive 
  adaptive refinement algorithms for the Poisson-Boltzmann equation. 
  Empirical results for the solvation free energy linear functional
  demonstrate that goal-oriented indicators are not sufficient on 
  their own to lead to a superior refinement algorithm.
  To remedy this, we propose a problem-specific marking strategy using
  the solvation free energy computed from the solution of the linear
  regularized Poisson-Boltzmann equation.
  The convergence of the solvation free energy using this 
  marking strategy, combined with goal-oriented refinement, compares
  favorably to adaptive methods using an energy-based error indicator.
  Due to the use of adaptive mesh refinement, it is critical to use 
  multilevel preconditioning in order to maintain optimal computational
  complexity.  We use variants of the classical multigrid method, which
  can be viewed as generalizations of the hierarchical basis multigrid and 
  Bramble-Pasciak-Xu (BPX) preconditioners.
\end{abstract}

\maketitle

\vspace*{-1.0cm}
\tableofcontents

\section{Introduction}
\label{sec:intro}

The Poisson-Boltzmann equation (PBE) 
is a widely used model for electrostatic interactions 
of charged bodies in dielectric media, such as molecules, ions, and colloids, 
and thus is of importance in many areas of science and engineering,
including biochemistry, biophysics, and medicine.
(See the classical texts~\cite{Mcqu73,Tanf61} for a derivation of the
PBE.) The importance of the PBE model is reflected by the popularity of 
software packages such as APBS~\cite{BSJH01}, CHARMM~\cite{CHARMM09},
DelPhi~\cite{DelPhi}, and UHBD~\cite{UHBD}, within the
molecular modeling communities.
It provides a high fidelity mean-field description of
electrostatic interactions and ionic distributions of solvated biomolecular 
systems in equilibrium. 
The partial differential equation itself is challenging to solve numerically
due to singularities of different orders at the positions of permanent point
charges and the presence of a dielectric interface. 

In this article, we develop an adaptive multilevel finite element
method for the PBE using {\em goal-oriented a posteriori error
indicators}.  This adaptive algorithm, which is a variant of that
studied for the PBE in \cite{HBW00, BHW00, BSJH01}, deviates substantially
from previous work in that the error indicator is based on a user
defined quantity of interest or goal. 
This is in contrast to traditional residual-based adaptive refinement algorithms (like
those developed for the PBE in~\cite{CHX07}) that drive-refinement to minimize the global
error measured in an energy-norm.
The goal-oriented refinement methodology has been successfully employed in 
a wide range of application areas, including fluids, elasticity, and fluid 
structure interaction~\cite{BaRa03}.
Despite these successes, we show that this methodology applied
directly to the PBE does not necessarily lead to a successful adaptive
algorithm. To remedy this issue we propose a novel marking strategy
which recovers the performance commonly seen in other applications.
This is the first time that this particular goal-oriented refinement strategy
has been applied to the PBE specifically, and molecular biophysics in general.

At the core of any adaptive finite element approach are the iterative
methods used to solve the discretized equation.  However, due to the
ill-conditioning of the linear systems arising from the discretization
of the PBE, the convergence rate of traditional iterative solvers is
significantly deteriorated.  To remedy this, we combine modern
Bramble-Pasciak-Xu (BPX)-type multilevel preconditioners with the
goal-oriented adaptive algorithm mentioned above.  When applied to the
PBE, our results demonstrate that the overall algorithm is accurate,
highly efficient and scalable with respect to the number of levels in
the adaptive hierarchy.

An outline of the article is as follows.
In \autoref{sec:pbe}, we give a brief overview of the Poisson-Boltzmann 
equation, and describe the most useful formulations for modeling and 
numerical simulation, such as the regularized formulations 
described in~\cite{CHX07,HYZZ10,CBO11}.
We also discuss the solvation free energy functional corresponding
to a given reaction potential, which will form the basis of our goal-oriented
error indicators developed later in the article.
We describe adaptive finite element methods in \autoref{sec:afem},
including weak formulation of the regularized PBE, 
discretization by finite element methods, 
and adaptive algorithms driven by {\em a posteriori} error indicators.
In \autoref{sec:goal}, we describe a particular class of error indicators
known as {\em goal-oriented} indicators, and describe 
several indicators designed for the PBE.
In \autoref{sec:prec}, we discuss a local multigrid algorithm
used to precondition an iterative Krylov method 
for solving the linear systems arising from adaptive mesh refinement.
With some care, these methods enable an algorithm whose complexity
is close to optimal.
The results from a sequence of numerical experiments using
the Finite Element ToolKit (\FETK{}) are presented in \autoref{sec:num}.
These results
highlight the efficacy of the goal-oriented error indicator for the
Poisson-Boltzmann problem, as well as the utility of
the linear solver strategy combined with the adaptive algorithm, driven
by the goal-oriented error indicator.
We draw some conclusions in \autoref{sec:conc}.

\section{The Poisson-Boltzmann Equation}
\label{sec:pbe}
\begin{figure}[t]
\begin{center} 
\mbox{\includegraphics[width=3.0in]{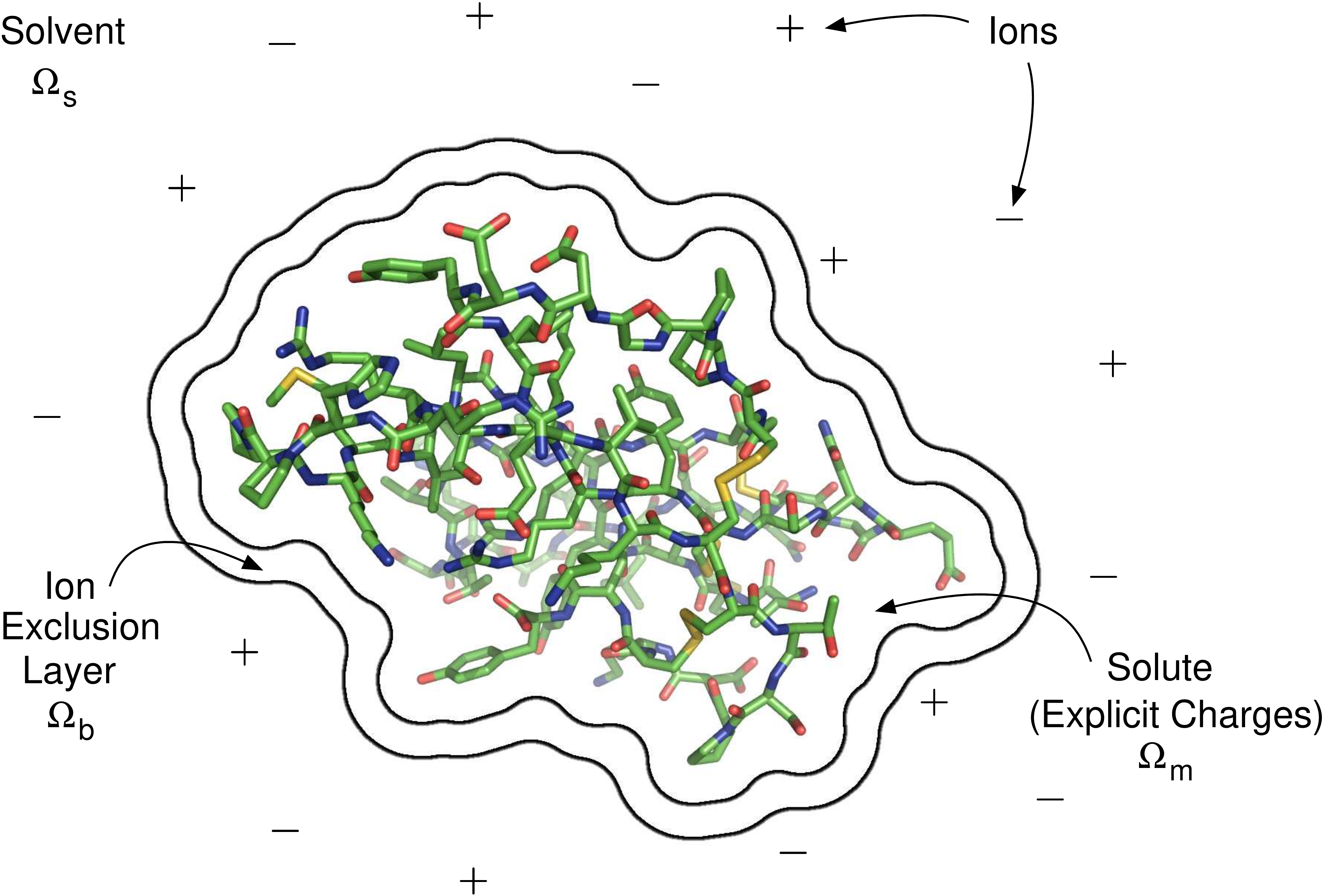}}
\end{center}
\caption{\label{fig:pbe-domain}
Debye-H\"uckel model of a charged biological structure immersed 
in a solvent containing mobile ions
}
\end{figure}

The Poisson-Boltzmann equation (PBE) is a second-order 
nonlinear partial differential equation whose solution gives the electrostatic
potential, $\phi(x)$, for a solute molecule immersed in an implicitly defined solvent.  Using a mean-field approximation,
the solvent is treated as a bulk medium where ions are distributed according to the
Boltzmann distribution. Figure~\ref{fig:pbe-domain} is a schematic representation of the domain of the PBE, denoted by $\Omega$.
The innermost region, $\Omega_m$, contains the explicitly represented solute molecule.  The outer
region, $\Omega_s$, is the bulk solvent and contains the implicit solvent ions.
Between $\Omega_s$ and $\Omega_m$ is the ion exclusion layer, which 
separates the solute from the solvent ions, and has a width dependent on
the size of the solvent ions.
For simplicity, we will assume the solvent ions are small and the ion exclusion layer 
can be neglected.  Hence, the interface between
the solute and solvent is a surface, denoted by 
$\Gamma = \bar{\Omega}_m \cap \bar{\Omega}_s$.  The shape of the surface is
governed by the short-range repulsive van der Waals interactions, which 
prevent the solvent from penetrating the solute. The precise definition 
of the surface varies depending on the model~\cite{BBC06}.

The PBE for a 1:1 electrolyte (e.g., sodium chloride) is
\begin{equation}\label{eqn:npbe}
\begin{split}
-\nabla\cdot\epsilon(x) \nabla u(x) + \bar{\kappa}^2(x)\sinh\left(u(x)\right)
 & =  \frac{4 \pi e_c}{k_B T} \sum_{i=1}^P q_i \delta(x-x_i), \quad x \in \Omega_m \cup \Omega_s, \\
u(\infty) & = 0, \\
\jump{\epsilon(x)\frac{\partial u(x)}{\partial n}}
& = 0, \quad x \in \Gamma,
\end{split}
\end{equation}
where $u(x)=e_c \phi(x)/k_B T$ is the dimensionless potential, $e_c$ is the 
charge of an electron, $k_B$ is Boltzmann's constant, and $T$ is the 
temperature.  Here, $\sjump{\cdot}$ denotes the jump across the interface
\begin{equation}
\jump{f(x)}
= \lim_{\zeta\rightarrow 0} f(x+\zeta n) - f(x-\zeta n)
\end{equation}
and $n$ is the outward pointing normal of $\partial\Omega_m$.
The dielectric function $\epsilon(x)$ jumps one or two orders of magnitude
at the interface $\Gamma$. For example, commonly used values are
$\epsilon(\Omega_m) = \epsilon_m = 2$ 
and $\epsilon(\Omega_s) = \epsilon_s = 80$. The modified Debye-H\"uckel 
parameter, $\bar{\kappa}$, has a similar discontinuity, with
$\bar{\kappa}(\Omega_m) = 0$ and 
$\bar{\kappa}(\Omega_s) = \bar{\kappa}_s > 0$.  The fixed ions within the
solute are represented by a sum of Dirac delta distributions, with fixed 
charge centers, $x_i$, and charges $e_c q_i$. This charge distribution
induces singularities in the electrostatic potential and has, until recently, 
proved to be difficult to treat numerically. 

To address this issue the PBE is reformulated so that the singularities
are explicitly removed~\cite{GDLM93,ZPVKL96}. 
Following~\cite{CHX07,HYZZ10,CBO11},
this is accomplished by writing the potential as a sum
of a singular term $u_c$ and a nonsingular remainder $u_r$.
The singular term is the Coulomb potential
\begin{equation}
u_c(x) = \frac{e_c }{\epsilon_m k_B T}\sum_{i=1}^P \frac{q_i}{|x-x_i|},
\end{equation}
which satisfies the Poisson equation
\begin{equation}\label{eqn:coulomb}
\begin{split}
-\nabla\cdot\epsilon_m\nabla u_c(x) = 
 \frac{4 \pi e_c}{k_B T} \sum_{i=1}^P q_i \delta(x-x_i) \quad \mbox{for} \quad x \in \Omega, \\
u_c(\infty) = 0.
\end{split}
\end{equation}
There are numerous fast algorithms (with linear or near linear complexity)
for evaluating $u_c$ on a set of quadrature points, such as fast
multipole~\cite{GrRo87}, multilevel
summation~\cite{SkTe02, Hardy06} and particle mesh Ewald~\cite{P3ME}.

Substituting $u=u_r+u_c$ into the PBE (Eq.~\ref{eqn:npbe}) 
gives a modified form of the PBE, 
which was termed in~\cite{CHX07} as the Regularized PBE (RPBE):
\begin{equation}\label{eqn:rpbe}
\begin{split}
-\nabla\cdot\epsilon(x) \nabla u_r(x) + & \bar{\kappa}^2(x)\sinh\left(u_r(x) + u_c(x)\right) \\
 & =  \nabla\cdot(\epsilon(x)-\epsilon_m) \nabla u_c(x), \quad x\in\Omega_m\cup\Omega_s, \\
u_r(\infty) & = 0, \\
\jump{\epsilon(x)\frac{\partial u_r(x)}{\partial n}}
& = (\epsilon_m-\epsilon_s)\frac{\partial u_c(x)}{\partial n},
\quad x \in \Gamma. 
\end{split}
\end{equation}
Note that, because
both $\bar{\kappa}(x)$ and $(\epsilon(x)-\epsilon_m)$ are zero for
$x\in\Omega_m$ and the
centers of the atoms in the solute are well separated from $\Gamma$,
the singular function, $u_c$, is never evaluated near the singularities. 
This formulation was used in~\cite{CHX07} to develop continuous {\em a priori}
$L^{\infty}$ estimates of solutions, and subsequently to show existence and 
uniqueness of solutions of the PBE.
Furthermore, the authors established discrete {\em a priori} $L^{\infty}$
estimates for Galerkin solutions, making possible quasi-optimal {\em a priori}
error estimates, as well as a provably convergent adaptive finite
element method for the RPBE.

Recently, an alternative 3-term splitting of the PBE has been proposed
in~\cite{HYZZ10} which addresses the inherit subtractive cancellation
in the reconstruction of the electrostatic potential used by the RPBE.
In this article, the authors establish mathematical results for the 
alternative splitting, including continuous and discrete
{\em a priori} $L^{\infty}$ estimates, existence and uniqueness of solutions,
quasi-optimal {\em a priori} error estimates, and a convergent adaptive
finite element method (AFEM).
(Whereas in~\cite{CHX07} only AFEM convergence was shown, it was shown 
in~\cite{HYZZ10} that AFEM is a contraction for the RPBE, using a new
AFEM convergence framework for nonlinear problems developed in~\cite{HNT09}.)
The 3-term splitting decomposes the electrostatic potential into 
\begin{equation}
u(x) =
\left\{
\begin{array}{ll}
u_3(x) + u_c(x) + u_h(x) & \mbox{in}\; \Omega_m \\
u_3(x)                   & \mbox{in}\; \Omega_s
\end{array}
\right.,
\label{eqn:3-term_decomp}
\end{equation}
where $u_c$ is the Coulomb potential, but here it is restricted to the
subdomain $\Omega_m$. The harmonic term, $u_h$, is defined as the solution to 
\begin{alignat}{2}
-\nabla^2 u_h(x) & = 0 &\quad& \text{in $\Omega_m$} \label{eqn:r3pbe-harm}\\
u_h(x) & = -u_c(x)     &\quad& \text{on $\Gamma$}.
\label{eqn:3-term_harm}
\end{alignat}
Applying the definitions of $u_c$, $u_h$ and substituting into the
PBE (Eq.~\ref{eqn:npbe}), one obtains an equation for $u_3$,
\begin{alignat}{2}
-\nabla\cdot(\epsilon(x) \nabla u_3(x)) + \bar{\kappa}^2(x) \sinh(u_3(x)) = 0 &\quad& \text{in $\Omega$}, \\
\jump{\epsilon(x)\frac{\partial u_3(x)}{\partial n}}
= \epsilon_m \frac{\partial(u_c(x)+u_h(x))}{\partial n}  &\quad&\text{on $\Gamma$}, \\
u_3(\infty) = 0. &\quad&  
\label{eqn:3-term_reg}
\end{alignat}
In contrast to the RPBE, this formulation avoids the subtractive 
cancellation since $u_3 = u$ in $\Omega_s$, and the Coulomb term is not
used to reconstruct the potential in the solvent subdomain.

For systems which are not highly charged, the variation in the potential
is relatively small, and the hyperbolic sine term is well approximated by
its linearization.  This approximation, which replaces $\sinh(u)$ with $u$, 
results in what is known as the linear Poisson-Boltzmann equation (LPBE), or 
the linear regularized Poisson-Boltzmann equation (LRPBE) in case of the 
RPBE.  Although this approximation reduces the ionic response of the 
solvent~\cite{BBC06,FZEV99}, it can significantly reduce the complexity of 
many numerical algorithms (e.g., boundary element and boundary integral 
methods)~\cite{LZHB07,BoFe04,LZHM08}.

One important use for the solution to the PBE is in the calculation of solvation
free energies. This quantity measures the thermodynamic work of moving the
solute molecule from a vacuum to a solvent environment.
The solvation free energy can be written as a sum of nonpolar and polar
contributions.  The nonpolar term depends on the solvent accessible surface 
area, excluded volume, and nonpolar forces which are typically assumed to
be independent of the electrostatic potential~\cite{IBR98,WaBa06}.
The polar term, $S$, is a
linear functional of the solution to the RPBE, $u_r$, (also known as
the reaction potential)~\cite{IBR98} and can be expressed
as
\begin{equation}\label{eqn:solvfe}
S(u_r) = \frac{1}{2} \int_\Omega u_r(x) \sum_{i=1}^P e_c \, q_i \delta(x-x_i) \;dx.
\end{equation}
For the 3-term splitting, the reaction potential is the sum of $u_3$ and $u_h$.

\section{Adaptive Finite Element Methods}
\label{sec:afem}
The finite element approach provides a natural framework for dealing with the
complex molecular surfaces which arise in the PBE.  Although there are modified 
finite difference methods which address this difficulty~\cite{YGW07}, finite 
element methods provide an attractive alternative when paired with an 
adaptive unstructured mesh designed to conform to the shape of the
solute molecule~\cite{LZHM08}. In 
this section, we present a general adaptive finite element
method for the regularized PBE, including the weak formulation, discretization,
solution using an inexact global Newton iteration, and adaptive refinement
procedure.  For more details on the finite element method,
see~\cite{StFi73,BrSc02,Brae01,ErGu04}.

\subsection{Weak Forms}
To give a well-defined weak formulation, the nonlinearity involving
exponentials must be controlled; in~\cite{CHX07,HYZZ10}, {\em a priori}
$L^{\infty}$ estimates are obtained for any solution to the RPBE, giving
almost everywhere pointwise bounds of the form: $\alpha \leqs u_r \leqs \beta$.
This leads to working with a well-defined solution space that consists 
of a non-empty, topologically closed, convex subset of $H^1(\Omega)$:
\begin{equation}
M_e := \{~ v \in H^1(\Omega) ~:~ \alpha \leqs v \leqs \beta
\mbox{ a.e. in } \Omega, v = u-u_c \mbox{ on } \partial\Omega ~\}.
\end{equation}
It is shown in~\cite{CHX07,HYZZ10} that there exists a unique solution 
to either regularized form of the RPBE in $M_e \subset H^1(\Omega)$. 
The weak formulation is: Find $u_r \in M_e$ such that
\begin{equation}
a(u_r,v)+b(u_r+u_c,v) = L(v) \quad\forall v \in H_0^1(\Omega)
\label{eqn:rpbe-weak}
\end{equation}
where
\begin{align}
a(u,v) & = (\epsilon\nabla u,\nabla v) \\
b(u,v) & = (\bar{\kappa}^2 \sinh(u),v) \\
L(v) & = (-(\epsilon-\epsilon_m)\nabla u_c,\nabla v).
\label{eqn:weak-form-components}
\end{align}
The linear functional $L(\cdot)$ is defined by integrating the right hand side of
Eq.~\ref{eqn:rpbe} by parts and applying the jump condition to eliminate the interface terms.

The weak form for the $3$-term split regularized PBE requires solving two
problems: first for the harmonic term on $\Omega_m$ and second for the 
split potential on the whole domain $\Omega$. Define the solution space to the harmonic
problem as
$M_h := \{ v \in H^1(\Omega_m) : v(x) = -u_c(x) \;\forall x \in \partial\Omega_m \}$.
Then the weak form of Eqns.~\ref{eqn:r3pbe-harm}-\ref{eqn:3-term_reg}
is: Find $(u_h,u_3) \in M_h \times M_e$ such that
\begin{equation}
a(u_3,v) + b(u_3,v) + a_m(u_h,w) = \langle g(u_h),v\rangle
\quad \forall (w,v) \in H_0^1(\Omega_m) \times H_0^1(\Omega)
\label{eqn:rpbe3-weak}
\end{equation}
where $a_m(\cdot,\cdot)$ is the restriction of the bilinear form to the
$\Omega_m$ subdomain and 
\begin{equation}
\langle g(u_h),v \rangle
= \int_{\partial\Omega_m} \epsilon_m \frac{\partial(u_c+u_h)}{\partial n} v 
  d x.
\end{equation}

\subsection{Solving}
Due to the hyperbolic sine, the RPBE has a strong nonlinearity. The discretized nonlinear
problems defined in Eqns.~\ref{eqn:rpbe-weak} and~\ref{eqn:rpbe3-weak} can be solved
using an inexact-Newton method~\cite{EiWa94}.
For brevity, we give details for Eq.~\ref{eqn:rpbe-weak}.
Define the weak residual functional to be
\begin{equation}
\langle R(u_r^h), v \rangle = L(v)-\left(a(u_r^h,v) + b(u_r^h+u_c^h,v)\right).
\label{eqn:weak-residual-rpbe}
\end{equation}
Here $u_r^h$ is the discrete solution satisfying the system of nonlinear equations
\begin{equation}
\langle R(u_r^h), v \rangle = 0 \quad \forall v \in V^h
\label{eqn:nonlinear-system-rpbe}
\end{equation}
where $V^h$ is the space of piecewise linear functions defined by the tetrahedral mesh.
Linearizing Eq.~\ref{eqn:nonlinear-system-rpbe} around $u_r^h$ results in
\begin{align}
J w^h & := \langle DR(u_r^h)w^h,v\rangle \nonumber \\
& = \frac{d}{d\epsilon}\left(\langle R(u_r^h + \epsilon w^h),v\rangle\right|_{\epsilon=0} \nonumber \\
& = -a(w^h,v) - b'(u_r^h+u_c;w^h,v) \quad \forall v \in V^h.
\label{eqn:rpbe-dr-define}
\end{align}
In the linear RPBE, $\sinh(u)$ is replaced by $u$, and  
$b'(u,v)=(\bar{\kappa}^2 u,v)$.
Newton's method defines the nonlinear update vector $s^h$ as the solution to
\begin{equation}
J s^h  = -R(u_r^h).
\label{eqn:newton-method}
\end{equation}
Given an initial guess $u_r^h \approx u_0$, the updated solution is defined as $u_1 = u_0 + s_h$.
This process can be repeated until a desired level of convergence is achieved.
An inexact-Newton method uses an iterative solve to find an approximate 
solution to Eq.~\ref{eqn:newton-method}, with a relatively large tolerance for 
the linear solve when far from the nonlinear solution.
However, as the exact solution to Eq.~\ref{eqn:nonlinear-system-rpbe} is 
approached, the linear solver tolerance is tightened so that quadratic 
convergence is achieved.

The computational complexity of the Newton solver is dominated by the
method used to solve the $N$ linear algebraic equations
\cite{BaRo82,Hack85} within each iteration.  
Multilevel methods provide an advantage in that they are
provably optimal or nearly optimal methods for solving these
systems~\cite{BaDu81,Hack85,Xu92a}.
The presence of geometrically complex discontinuities in the 
dielectric $\epsilon$ and
in the Debye-H\"uckel parameter $\bar{\kappa}$ in the PBE destroy
classical multilevel method efficiency, and can even cause divergence.
This is analyzed at length for the PBE 
in~\cite{Hols94d,HoSa93a,Hols94e}, 
and various techniques based on coefficient averaging and algebraic
enforcement of variational (Galerkin) conditions are examined.
Algebraic multilevel methods have been used successfully for many similar
problems; cf.~\cite{BaXu94,BaXu96,Bran86,CSZ94,CGZ97,BMR84,RuSt87,VMB94,VMB95}.
A fully unstructured algebraic multilevel approach is taken in
\FETK{}, more details are provided in Section~\ref{sec:prec}.

Starting from an initial mesh $\calT_0$, the adaptive mesh refinement procedure
builds a sequence of conforming meshes
$\calT_0, \calT_1, \ldots , \calT_l$~\cite{CHX07,HYZZ10,Dor96}.
This procedure is divided into four steps: \semr{SOLVE}, \semr{ESTIMATE},
\semr{MARK}, and \semr{REFINE}. In the \semr{SOLVE} step, a solution is computed on the current
mesh. Using this result, the \semr{ESTIMATE} step
computes elementwise error indicators and an estimate of the global error.
In a production environment, the procedure
terminates if the global error estimate is below some prescribed tolerance.  Otherwise,
the \semr{MARK} step selects elements for refinement. This step is crucial to the convergence
of the method.  
Finally, \semr{REFINE} subdivides the marked elements possibly subdividing additional unmarked
elements in order to produce a conforming mesh.
The refinement technique used in this article is longest edge bisection~\cite{Riva84a}.

\section{Error Indicators}
\label{sec:goal}
In this section we present \emph{a posteriori} error indicators for use 
in adaptive refinement. These estimators
are typically developed by considering the residual of the weak form.
For example, given
a finite element solution $u_r^h \in V^h$, the weak residual for the
linear RPBE is (compare to Eq.~\ref{eqn:weak-residual-rpbe})
\begin{equation}
\langle R(u_r^h), v \rangle  = L(v) - (\bar{\kappa}^2 u_c,v)
-\left(a(u_r^h,v) + (\bar{\kappa}^2 u_r^h,v)\right)
\label{eqn:weak-residual-lrpbe}
\end{equation}
for a given $v\in V$. For the remainder of this article, we restrict our
attention
to two classes of error indicators: energy-based and goal-oriented.
The first class estimates the error in the energy norm, although this idea can
be generalized to other norms, (e.g., the $H^1$ norm).  The second class, 
called goal-oriented, focuses on estimating the error in a user
specified quantity of interest or goal functional.  In the following sections,
we derive error indicators from both classes for the linear RPBE. 
For the goal-oriented indicator, the solvation free energy is used as
the target functional.

\subsection{Energy Norm Indicators}
A standard \emph{a posteriori} error indicator is based on bounding the error in the energy norm.
It is easily derived by breaking the weak residual into its elementwise components and integrating
by parts over each element~\cite{AiOd00}. This technique was used in~\cite{Cyr08}
to derive the following estimator for the linear RPBE
\begin{equation}\label{eqn:ener-indicator}
\eta_K^2(u^h) = h_K^2\|r_K\|_{L^2(K)}^2 + \frac14 h_{\partial K} \|r_{\partial K}\|_{L^2(\partial K)}^2,
\end{equation}
where
\begin{equation}
\begin{split}
r_K(x) = & (\nabla \cdot(\epsilon(x)-\epsilon_m) \nabla u_c(x) - \bar{\kappa}^2(x) u_c(x)) \\
& - (-\nabla \cdot \epsilon(x)\nabla u_r^h(x) + \bar{\kappa}^2(x) u_r^h(x)) \quad \quad \forall K \in \calT
\end{split}
\end{equation}
and
\begin{equation}
r_{\partial K}(x) = n_K \cdot \jump{(\epsilon(x)-\epsilon_m) \nabla u_c(x) + \epsilon(x) 
\nabla u_r^h(x)}_{n_K}
\quad \forall K \in \calT.
\end{equation}
This indicator gives a bound on the error measured in the energy norm ($\energy{v}^2 := a(v,v) +(\bar{\kappa}^2 v,v)$)
\begin{equation}
\energy{u_r-u_r^h}^2 \leq \sum_{K\in\calT} \eta_K^2(u_r^h).
\end{equation}
Bounding the error in other norms is possible. For example, 
in~\cite{CHX07,HYZZ10} a similar \emph{a posteriori} error indicator for the RPBE 
was shown to bound the error measured in the $H^1$ norm. Other efforts have focused
on formulating the RPBE as a first-order system least squares (FOSLS) problem, which has a natural
error estimate~\cite{BCCO10,CBO11b}.

\subsection{Goal-Oriented Indicators}

Key to the development of goal-oriented error indicators is relating the weak residual to the
error in the goal functional. For symmetric linear problems, a direct application of the
Riesz representation 
theorem shows that there exists a dual function that when paired with the weak residual gives the error 
in the goal~\cite{BeRa01}. The challenge is to approximate this function and utilize that approximation to develop error
indicators. However, for nonlinear problems, like Eq.~\ref{eqn:rpbe}, the definition of the dual function is
not so clear. The first part of this section discusses a strategy for defining the dual function for 
both the nonlinear RPBE and the three term splitting. 
Using this definition of the dual, two goal-oriented error indicators are proposed for the linear RPBE utilizing
the solvation free energy as the quantity of interest.

Following~\cite{BeRa01}, a dual function for the nonlinear RPBE can be defined by considering
the constrained minimization problem
\begin{equation}
u_r = \underset{u_r^*\in M_e}{\arg\min} \; S(u_r^*) \quad \mbox{subject to} \quad a(u_r,v)+b(u_r+u_c,v) = L(v)
\;\;\forall v \in H_0^1(\Omega),
\end{equation}
where $a(\cdot,\cdot)$, $b(\cdot,\cdot)$ and $L(\cdot)$ are specified in Eq.~\ref{eqn:weak-form-components}, and $S$ is the goal functional.
Note that because the RPBE constraint determines the solution uniquely, the minimization problem has the same solution. However, specifying the
minimization provides an additional mathematical framework to define the dual function. To see this, consider the Lagrangian associated
with the minimization problem
\begin{equation}
\Theta(u_r,w) = S(u_r) + \bigl(L(w) - \left(a(u_r,w)+b(u_r+u_c,w)\right)\bigr),
\end{equation}
where the Lagrange multiplier, $w \in H_0^1(\Omega)$, is also the dual function.
Taking the first variation of $\Theta$ with respect to $u$ gives the dual problem:
\begin{equation}\label{eqn:rpbe-dual-problem}
\mbox{Find} \,\, w \in H_0^1(\Omega), \,\, \mbox{such that} \,\, \langle D R(u_r)v,w\rangle = -S(v), \quad \forall v \in H_0^1(\Omega),
\end{equation}
where $\langle D R(\cdot)\cdot,\cdot\rangle$ was defined in Eq.~\ref{eqn:rpbe-dr-define}.
As discussed above, if $b(\cdot,\cdot)$ is linear in the first argument then Eq.~\ref{eqn:rpbe-dual-problem} simplifies 
to $a(v,w)+(\bar{\kappa}^2 v,w) = S(v)\;\;\forall v \in H_0^1(\Omega)$. For the linear problem the error in a goal functional $S(\cdot)$,
like the solvation free energy, is simply expressed in terms of the weak residual
\begin{eqnarray}
S(u_r-u_r^h) & = & a(u_r-u_r^h,w)+(\bar{\kappa}^2 (u_r-u_r^h),w) \nonumber \\
& =&  L(w)-a(u_r^h,w)-(\bar{\kappa}^2 (u_r^h+u_c),w).\label{eqn:solv-dual-err}
\end{eqnarray}
Thus, if $w$ is known, the error in $S(\cdot)$ is easily calculated. In the nonlinear case the
error in the goal satisfies
\begin{equation}
S(u_r-u_r^h) = L(w)-a(u_r^h,w)-b'(u_r^h+u_c;u_r^h,w)+ E,
\end{equation}
where $E$ is quadratic in the error in $u_r^h$~\cite{BeRa01}. 

For the three term splitting of the PBE, again we setup a constrained minimization problem to define the
dual. 
Using the notation from 
the weak form in Eq.~\ref{eqn:rpbe3-weak}, the corresponding Lagrangian is
\begin{multline}
\Theta_{3\mbox{-term}}(u_3,u_h;w_3,w_h) = S(u_3) + S_m(u_h) + \\ 
\langle g(u_h),w_3 \rangle - \left(a(u_3,w_3) + b(u_3,w_3) + a_m(u_h,w_h)\right)
\end{multline}
where $w_3 \in H_0^1(\Omega)$ and $w_h \in H_0^1(\Omega_m)$ are dual functions. The functional
$S_m(\cdot)$ is a restriction of the
original goal functional to the $\Omega_m$ domain.
Taking the first variation of $\Theta_{3\mbox{-term}}$ with respect to $u_3$ and $u_h$ gives
the dual problem
\begin{equation}
\begin{aligned}
a(v,w_3)+b'(u_3;v,w_3) & = S(v) & \forall v \in H_0^1(\Omega), \\
a_m(v,w_h) & = S_m(v)+\langle g'(u_h;v), w_3\rangle & \forall v \in H_0^1(\Omega_m)
\end{aligned}
\end{equation}
where
\begin{equation}
\langle g'(u_h;v),w_3 \rangle
= \int_{\partial\Omega_m} \epsilon_m \frac{\partial v}{\partial n} w_3 
  d x.
\end{equation}

\subsubsection{Goal-Oriented Error Indicators for the Linear RPBE}\label{sec:lpbe-goal}
To make the application of the dual functions in error indicators more concrete, we present
two goal-oriented indicators for the linearized RPBE. As an example we will focus on accurate computation of
the solvation free energy (see Eq.~\ref{eqn:solvfe}).  Unfortunately, $S(\cdot)$ is not
bounded on $H_0^1(\Omega)$ due to the inclusion of delta distributions. A common approach
to circumventing this issue is to use a \emph{mollified} version of the
functional~\cite{AiOd00,BaRa03}. In this case, the mollified solvation free
energy is
\begin{equation}
S(u_r) \approx S^\sigma(u_r) = \frac{1}{2}\int u_r(x) \sum_{i=1}^P e_c \, q_i \theta(|x-x_i|,\sigma) \diff x, 
\end{equation}
where $\theta$ is a locally supported function defined such that
\begin{equation}
\underset{\sigma \rightarrow 0}{\lim} \int \theta(|x|,\sigma) f(x) \diff x = \int \delta(x)f(x) \diff x= f(0).
\end{equation}
One possible choice for $\theta$ is the step function
\begin{equation}\label{eqn:moll-kernel}
\theta(r,\sigma) = \left\{
\begin{array}{cc}
B_\sigma^{-1} & r \leq \sigma, \\
0 & r > \sigma,
\end{array}
\right.
\end{equation}
where $B_\sigma$ is the volume of a ball of radius $\sigma$.

A simple error indicator suggested by Eq.~\ref{eqn:solv-dual-err} is to first solve the dual problem using
the same approximation space as the primal.
This approximate dual could then be substituted for $w$ in Eq.~\ref{eqn:solv-dual-err}
to compute the value of the indicator.  However, if the same finite dimensional space is used for solving both the dual
problem and the primal problem, then, because of Galerkin orthogonality, Eq.~\ref{eqn:solv-dual-err}
will be zero (see~\cite{BaRa03}). A remedy is to instead solve the dual problem using a finer approximation space,
$U^h \subset H_0^1(\Omega)$.  One convenient choice is to maintain the same mesh and use higher
order polynomials for $U^h$. In the examples below, $V^h$
is the space of piecewise linear polynomials and $U^h$ is the space of piecewise quadratics.
Let the finer resolution solution of the dual problem be denoted $w^{h,2}$. Substituting $w \approx w^{h,2}$
into Eq.~\ref{eqn:solv-dual-err} yields
\begin{multline}
S^\sigma(u_r-u_r^h) \approx  L(w^{h,2}-P_{h,2}^h w^{h,2}) \\
- \left(a(u_r^h,w^{h,2}-P_{h,2}^h w^{h,2})
 + (\bar{\kappa}^2 (u_r^h+u_c),w^{h,2}-P_{h,2}^h w^{h,2})\right)
\label{eqn:solv-err-eqn}
\end{multline}
Where $P_{h,2}^h$ is a convenient projection (e.g., nodal injection) of the fine space $U^h$
onto $V^h$. The choice of the projection operator will 
affect the quality of the indicator. Decomposing the error into its elementwise
contributions gives
\begin{multline}\label{eqn:goal-quad-ind-first}
S^\sigma(u_r-u_r^h) \approx \sum_K L_K(w^{h,2}-P_{h,2}^h w^{h,2}) 
- a_K(u^h,w^{h,2}-P_{h,2}^h w^{h,2})\\ -(\bar{\kappa}^2 (u_r^h+u_c),w^{h,2}-P_{h,2}^h w^{h,2})_K 
\leq \sum_K \eta_K(u^h,w^{h,2}),
\end{multline}
where the subscript $K$ indicates the restrictions of the linear functional, bilinear functional or
inner product to element $K$ and
\begin{equation}\label{eqn:goal-quad-ind}
\begin{split}
\eta_K(u^h,w^{h,2}) =& \int_K \biggl| -(\epsilon(x)-\epsilon_m) \nabla u_c(x) \cdot \nabla(w^{h,2}(x)-w^h(x)) \\
& - \bar{\kappa}^2(x) u_c(x)(w^{h,2}(x)-w^h(x)) \\
& -\epsilon(x) \nabla u_r^h (x) \cdot \nabla(w^{h,2}(x)-w^h(x)) \\
& - \bar{\kappa}^2(x) u_r^h(x)(w^{h,2}(x)-w^h(x))\biggl| \diff x.
\end{split}
\end{equation}
Here, the absolute value of the integrand in Eq.~\ref{eqn:goal-quad-ind-first} 
has been taken over each element.  In the numerical experiments in 
section~\ref{sec:num}, Eq.~\ref{eqn:goal-quad-ind} is referred to as the
the ``goal-quadratic'' error estimator used by the various adaptive refinement 
marking strategies.

The error indicator in Eq.~\ref{eqn:goal-quad-ind} requires solving a dual problem which is substantially
larger than the primal problem. To alleviate this issue, we develop a second goal-oriented error estimator
that finds an approximation to the dual in $V^h$ (which is the same space as the primal problem). 
The error in the goal is estimated by solving many local elementwise boundary value problems.
The technique proposed here is similar to the development of the \emph{equilibrated residual method}
for computing goal-oriented estimators~\cite{OdPr01,PrOd99}.
However, the less accurate but simpler \emph{element residual method} (ERM), as discussed in~\cite{AiOd00},
is used.  Using the parallelogram law, the error in the linear functional can
be rewritten~\cite{AiOd00,OdPr01,PrOd99} as
\begin{equation}\label{eqn:parallel}
\begin{split}
S^\sigma(u_r-u_r^h) & = a(u_r-u^h_r,w-w^h) + (\bar{\kappa}^2 (u_r-u_r^h),w-w^h) \\
                   & = \frac14\energy{(u_r-u_r^h) + (w-w^h)}^2- \frac14\energy{(u_r-u_r^h) - (w-w^h)}^2
\end{split}
\end{equation}
where again $\energy{v}^2 = a(v,v) + (\bar{\kappa}^2 v,v)$ is the square of
the energy norm.
Define elementwise error functions $\phi_K = u_r-u_r^h|_K$ and $\psi_K = w-w^h|_K$ to be computed by the element residual method.
The error in the primal problem on element $K$ is approximated by the solution to
\begin{equation}
a_K(\phi_K,v)+(\bar{\kappa}^2 \phi_K, v)_K = \langle R(u_r^h), v \rangle_K
   + \int_{\partial_K} f_K^u v(x) \diff s  \quad \forall v \in H^1(K),
\end{equation}
where 
\begin{equation}
\langle R(u_r^h), v \rangle_K = L_K(v) - a_K(u_r^h,v)-(\bar{\kappa}^2(u_r^h+u_c),v)_K,
\end{equation}
and
\begin{equation}\label{eqn:erm-fluxes}
f_K^u = \bigl(\epsilon(x)\nabla u_r^h(x)+(\epsilon(x)-\epsilon_m) \nabla u_c(x)\bigr) \cdot n_K.
\end{equation}
Similarly, the error in the dual problem on element $K$ is approximated by the solution to
\begin{equation}
\begin{split}
a_K(v,\psi_K)+(\bar{\kappa}^2 v, \psi_K)_K = S_K^\sigma(v) - a_K(v,w^h) -(\bar{\kappa}^2 v,w^h)_K \\
+ \int_{\partial_K} f^w_K v(x) \diff s  \quad \forall v \in H^1(K),
\end{split}
\end{equation}
where
\begin{equation}\label{eqn:erm-fluxes2}
f_K^w = \epsilon(x)\nabla w^h(x) \cdot n_K.
\end{equation}
For a derivation of these equations see~\cite{Cyr08}. Explicitly
stated, the ``goal-linear'' error estimator used in the numerical experiments 
in section~\ref{sec:num} is given by
\begin{equation}\label{eqn:goal-linear-indicator}
\eta_K(u_r^h,w^h) = \frac14\energy{\phi_K + \psi_K}_K^2 - \frac14\energy{\phi_K-\psi_K}_K^2.
\end{equation}

\section{Multilevel Preconditioning}
\label{sec:prec}
As the mesh is refined, the conditioning of the linear system
deteriorates, and preconditioning is necessary to accelerate
the convergence of iterative solvers (e.g., the conjugate gradient 
method). The challenge in designing an efficient preconditioner is 
balancing the cost of applying the preconditioner with its effectiveness
in improving the conditioning of the underlying system.

In a given finite element mesh, at level $j$, we denote the set of
nodes and its cardinality by $\mathcal{N}_j$ and $N_j$,
respectively. We call the set of nodes introduced precisely at level
$j$ the \emph{fine} nodes, and denote them $\mathcal{N}_j^f$. As the mesh
is refined, $\mathcal{N}_j^f$ is appended to $\mathcal{N}_{j-1}$,
leading to the following hierarchy
of nodes:
\begin{equation*}
\mathcal{N}_j = \mathcal{N}_{j-1} \bigcup \mathcal{N}_j^f,\quad j=1, \ldots, J,
\end{equation*}
where $N_J = N$. 

In the local mesh refinement setting, the way the coarse and fine nodes
are processed plays a central role in determining the overall efficiency
of a preconditioner.  If the computational cost
per level can be maintained proportional to $N_j - N_{j-1}$, or
slightly larger, then total cost will be order $N$, and the resulting
preconditioner is said to have optimal
computational complexity per iteration.  If the resulting preconditioned
system has a bounded condition number (independent of problem size), a solution
can be obtained using an iterative method with a bounded number of iterations.
Hence, the combination of optimal per iteration complexity with a bounded
condition number leads to a solver with optimal overall complexity.

In this article, we restrict the presentation of local multilevel
preconditioning to a purely geometric (node based) perspective because
the computational complexity is exactly governed by the number of
nodes processed by the preconditioner at each level.  The local
multilevel preconditioners of interest can be classified into two groups: 
multiplicative local multigrid (MG)~\cite{Bast93,BaWi94,BPWX91,Riva84b}
and additive local MG~\cite{BPWX91,BPX90}.  The additive local MG
preconditioner is often called the Bramble-Pasciak-Xu (BPX)
preconditioner in the literature.  In this article, we report on only
the multiplicative variants. We use the
term ``classical'' to refer to the application of a preconditioner
in the uniform refinement setting.  There is abundant literature on MG
preconditioners~\cite{BHM00,TOS01,Xu92a,Yavn06} and local MG 
(e.g., see~\cite{BHW98,Witt04} for a review).

Proofs demonstrating the optimality of classical MG and classical 
additive MG preconditioners rely on a geometric increase 
in the number of nodes per level. This is because the cost per iteration 
of these classical preconditioners is proportional to $N_j$ 
(not $N_j - N_{j-1}$ ) per level, which results in suboptimal complexity 
if $\sum_{j=1}^J c_j \, N_j$ is not~$\mathcal{O}(N)$.  This occurs 
frequently in the local refinement setting due to the slow increase in 
the number of nodes between levels.

The 
family of
hierarchical basis (HB) preconditioners, developed 
by Bank, Dupont, and Yserentant~\cite{Bank96,BDY88,Yser93},
maintain a per-level cost proportional to $N_j - N_{j-1}$, by 
only processing (smoothing) the fine nodes at each level. Although the cost 
per iteration is optimal, HB preconditioners do not achieve a uniformly bounded
condition number, and suffer from $\mathcal{O}(J^2)$ and $\mathcal{O}(2^J)$
iterations in two- and three-dimensions, respectively. To address this 
deficiency, we investigate local MG preconditioners which process a larger 
set of nodes, $\mathcal{X}_j$, but still maintain a cost which is proportional 
to $N_j - N_{j-1}$ at each level. Hence, we seek a set $\mathcal{X}_j$ such
that
\begin{equation*}
\mathcal{N}_j^f \subset \mathcal{X}_j \subset \mathcal{N}_j,
\end{equation*}
with cardinality, $X_j$, which is proportional to $N_j - N_{j-1}$.  At
the same time, $\mathcal{X}_j$, should be large enough that the resulting
system has a bounded condition number, leading to a solver with optimal 
overall complexity. Aksoylu and Holst~\cite{AkHo06} showed that this is 
possible even for three-dimensional local refinement routines.

In the local mesh refinement setting, Aksoylu, Bond, and
Holst~\cite{ABH03} studied the implementation and algebraic aspects
(e.g., matrix representations) of multilevel preconditioners.  
Subsequent articles provide a comprehensive overview of local MG
preconditioners with various emphases: for a theoretical treatment, see
\cite{ABH04TR,AkHo05TR1,AkHo05TR2}; for
optimality analysis in three-dimensional local refinement routines,
see \cite{AkHo06}; for surface mesh applications in computer graphics,
see \cite{AkKhSc05}.

\subsection{Local multigrid preconditioners}
As mentioned in the previous section, the fundamental difference between 
classical and local MG preconditioners is the smoothing operation.
In classical MG, the smoother acts on all degrees of freedom on 
every level.  In contrast, local MG only smooths a small subset,
typically a neighborhood, of the fine degrees of freedom. Pseudo-code for
a local MG V-cycle is provided in Algorithm \ref{alg:localMG}.
\begin{algorithm} \label{alg:localMG}
\renewcommand{\arraystretch}{1.0}
Local multigrid V-Cycle: \\
$u_{[j]} = \mbox{Vcycle}(u_{[j]}, f_{[j]}, j)$ \\
\begin{tabular}{lll}
0) & If $j = 1$, solve $A_{[j]} u_{[j]}= f_{[j]}$ & coarsest level solve\\
   & \quad and return $u_{[j]}$; & \\
1) & $u_{[j]} = \mbox{Smooth}(u_{[j]}, A_{[j]}, f_{[j]}, \mathcal{X}_j, s_1)$; & smooth $s_1$ times on $\mathcal{X}_j$ \\
2) & $r_{[j-1]} = I_{[j]}^{[j-1]} \left(f_{[j]} - A_{[j]} u_{[j]}\right)$; & restrict residual \\
3) & $e_{[j-1]} = 0$; & set $e_{[j-1]}$ to zero \\
4) & $e_{[j-1]} = \mbox{Vcycle}(e_{[j-1]}, r_{[j-1]}, j-1)$; & coarse-grid recursion \\
5) & $u_{[j]} = u_{[j]} + I_{[j-1]}^{[j]} e_{[j-1]}$; & add interpolated correction \\
6) & $u_{[j]} = \mbox{Smooth}(u_{[j]}, A_{[j]}, f_{[j]}, \mathcal{X}_j, s_2)$; & smooth $s_2$ times on $\mathcal{X}_j$ \\
7) & return $u_{[j]}$ & \\
\end{tabular}
\end{algorithm}
The set of nodes used by the underlying method's smoothing operation
defines the type of preconditioner. In this article, we introduce three 
different sets of nodes to be used in the local MG preconditioner. On a
given refinement level, the \emph{marked region} is the set of elements 
which have been marked for refinement. The \emph{refinement region} is
union of the newly introduced elements as a result of the refinement
and closure procedure. Figure \ref{fig:nodeTypes} depicts a
two-dimensional refinement region
formed by a local refinement routine which consists of quadrasection
closed by bisection (the so-called \emph{red-green} refinement). The 
sets of nodes used to define the preconditioners are as follows:
\begin{figure}[ht]
\begin{center}
\includegraphics[width=.75\textwidth]{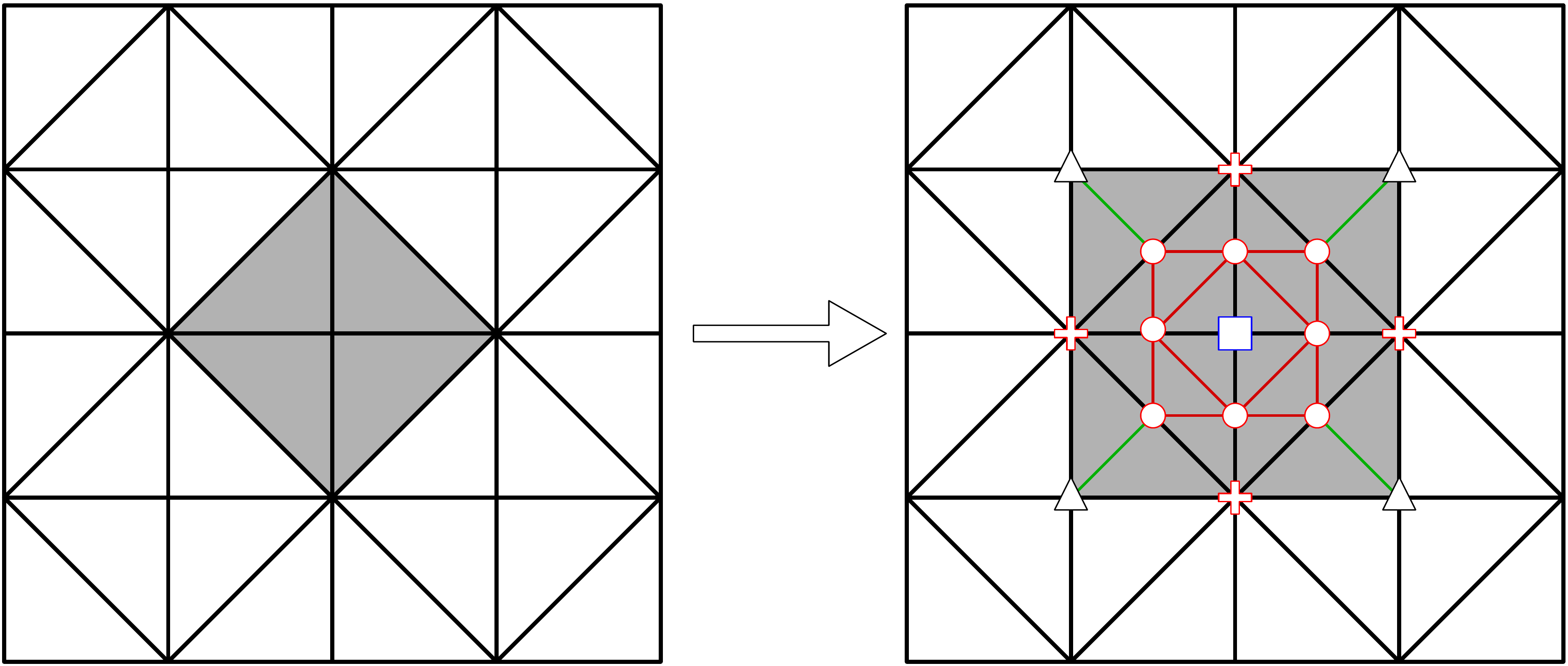}
\scriptsize
\begin{tabular}{lllll}
HB: & ${\color{red} \circ}$. &  &  \\
BPX:     & ${\color{red} \circ}$, & ${\color{blue} \square}$. &    & \\
BEK:     & ${\color{red} \circ}$, & ${\color{blue} \square}$, & ${\color{red} +}$. & \\
ONERING: & ${\color{red} \circ}$, & ${\color{blue} \square}$, & ${\color{red} +}$, & $\triangle$.
\end{tabular}
\end{center}
\caption{An example of red-green refinement 
(quadrasection-bisection) in two dimensions. On the left, the mesh before 
refinement, with the marked region shaded in gray. On the 
right, the mesh after refinement (red edges) and closure (green edges), with
the refinement region shaded in gray.  The labels 
indicate which nodes are assigned to each of the preconditioner smoothing 
sets, $\mathcal{X}_j$.
\label{fig:nodeTypes}}
\end{figure}
\begin{itemize}
\item $\mathcal{X}_j\text{-HB}$: The set of fine nodes on level $j$, i.e., 
$\mathcal{N}_j^f$~\cite{Bank96,BDY88,Yser93}.

\item $\mathcal{X}_j\text{-BPX}$: The set of nodes whose corresponding
basis functions have support entirely contained in the
\emph{refinement region}~\cite{BPWX91}.

\item $\mathcal{X}_j\text{-BEK}$: The set of nodes whose corresponding
  basis functions have non-empty intersection with the \emph{marked
    region}.  This set is named after the Bornemann-Erdmann-Kornhuber
  type refinement~\cite{BEK93} routine. It consists of fine
  nodes and their immediate neighboring coarse nodes in the
  \emph{marked region}.  For P1 elements, this set can be inferred
  from the nonzero pattern of the prolongation operator.

\item $\mathcal{X}_j\text{-ONERING}$: The set of nodes whose
  corresponding basis functions have non-empty intersection with the
  \emph{refinement region}~\cite{ABH03,BoYs93,BrPa93,DaKu92}. This set
  consists of fine nodes and their immediate neighboring coarse nodes
  in the \emph{refinement region}. For P1 elements, this set can be
  inferred from the nonzero pattern of the coarse-fine subblock of the
  stiffness matrix.
\end{itemize}

As an example, we have labeled the nodes in Figure \ref{fig:nodeTypes}
to show which nodes are in each of the sets described above. The set
corresponding to the classical multigrid preconditioner contains all of 
the nodes on each level.

We should note that the practical implementation of the various local 
MG preconditioners varies significantly depending on the particular
preconditioner. Special care must be taken in order to achieve optimal
computational as well as storage complexities. The implementation
aspects of how to construct optimal complexity preconditioners are
studied in more detail in~\cite{ABH03}.

\section{Numerical Experiments}
\label{sec:num}
\subsection{Adaptive Refinement}
In this section, the effectiveness of goal-oriented mesh refinement is compared
to refinement using the energy-based error indicator.  Of interest is computing
the solvation free energy of the $921$-atom Fasciculin-$1$
 protein~\cite{DMBF92},
using the solution of the linear RPBE.  All tests were performed
using \FETK{}~\cite{Hols01a}.

In order to solve the RPBE, a definition of the molecular surface and a mesh
conforming to that surface is needed. Various definitions of the molecular
surface have been proposed in the literature,
and the particular value obtained for the solvation free energy will depend
strongly on the surface geometry~\cite{BBC06}. However, the performance of the
algorithms proposed here are insensitive to the choice of molecular surface,
as long as
it is sufficiently smooth, and the underlying mesh is conforming.
Historically, generation of the
mesh conforming to the surface was a great impediment to 
using finite elements for solving the PBE, and only recently, with the 
development of tools like PDB2PQR and GAMer, has molecular meshing become a
routine task.  The first step is to prepare the structure using 
PDB2PQR~\cite{DNMB04}, which adds missing hydrogens, assigns charges, 
and specifies a radius for each atom in the protein.  Next, the resulting PQR
file is passed to GAMer~\cite{HYZZ10,HMYT09,YHCM08,YHM08}, which 
produces a tetrahedral mesh conforming to the shape of the protein.
Finally, this mesh is used by \FETK{} to solve the linear RPBE, and compute the
solvation free energy.

\subsubsection{Global Marking Strategy}
The 
adaptive refinement 
algorithm 
creates a sequence of grids $\calT_0,\calT_1,\ldots\calT_l\ldots$ based on an error indicator.
Critical to this algorithm is the third step MARK. The goal of this step is to select elements for refinement.
There are several different choices for marking strategies~\cite{AiOd00,BaRa03}.  
The strategy used here is to mark all elements in the $l^{th}$ refinement level that satisfy
\begin{equation}\label{eqn:pbe-simple-mark}
\eta_K > \gamma\; \underset{T \in \calT_l}{\max}\; \eta_T \quad \forall K \in \calT_l,
\end{equation}
where $\gamma \in (0,1)$. This criteria yields no refinement for $\gamma = 1$ and uniform refinement for $\gamma = 0$. 

\begin{figure}
\begin{center}
\includegraphics[width=.45\textwidth]{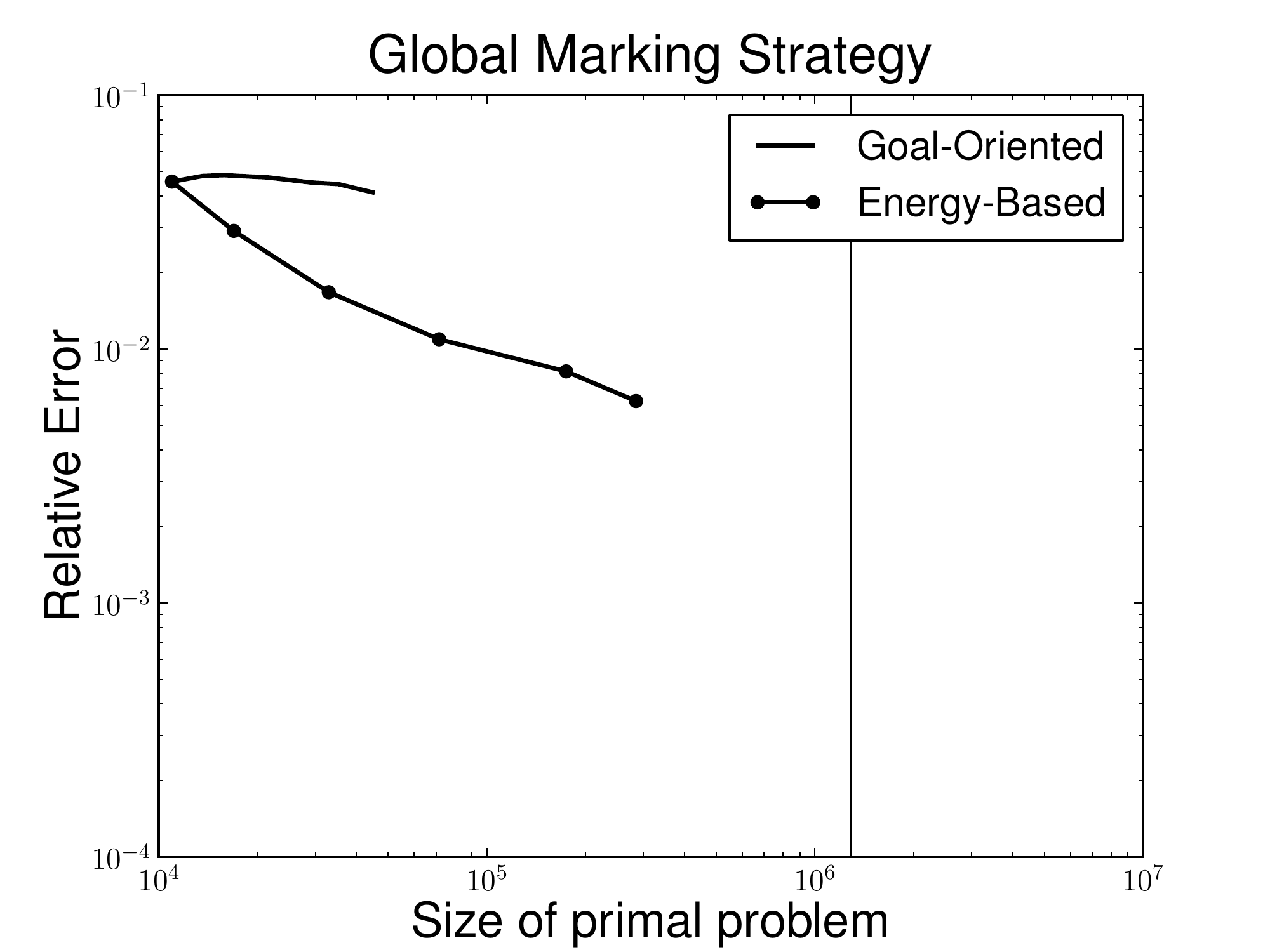}
\end{center}
\caption{Convergence of the solvation free energy for Fasciculin-$1$ with both goal-oriented and energy-based
indicators using the global marking strategy.}
\label{fig:global-sfe-1fas}
\end{figure}
 
The convergence of the solvation free energy for Fasciculin-1 using the global marking strategy with both
energy-based (Eq.~\eqref{eqn:ener-indicator}) and goal-oriented indicators (Eq.~\eqref{eqn:goal-quad-ind}) can be
seen in Fig.~\ref{fig:global-sfe-1fas}.
The figure shows the relative error in the solvation free energy as a function of the
number of unknowns in the primal problem. The error in the solvation free energy is estimated by
computing a high resolution solution to the PBE using uniform mesh refinement. The vertical
line near $1.5\times 10^6$ unknowns marks the size of this reference solution. Notice that this
figure differs from traditional finite element convergence plots which show error as a function
of element radius. Since the resolution of an adaptively refined mesh can greatly vary over the
spatial domain, the element radius is not an appropriate measure of the resolution of the mesh.
Furthermore, the order of accuracy of the approximation is based on the type of
basis functions used, and thus is the same for both uniform and adaptive
refinement. The benefit of adaptive refinement is that the mesh can be 
refined in regions that heavily contribute to the error, resulting in higher
accuracy with fewer total degrees of freedom.

In Fig.~\ref{fig:global-sfe-1fas}, the lower line shows the convergence of the solvation free energy for several
levels of adaptive refinement using the energy-based
indicator with the global marking strategy. This scheme makes steady progress
to the correct solvation free energy.  On the
other hand, the upper line shows the results using goal-oriented indicators with
the global marking strategy. This scheme performs very poorly.

\begin{figure}
\begin{center}
\begin{tabular}{cc}
  {Energy-Based} & {Goal-Oriented} \\
  \includegraphics[width=0.40\textwidth]{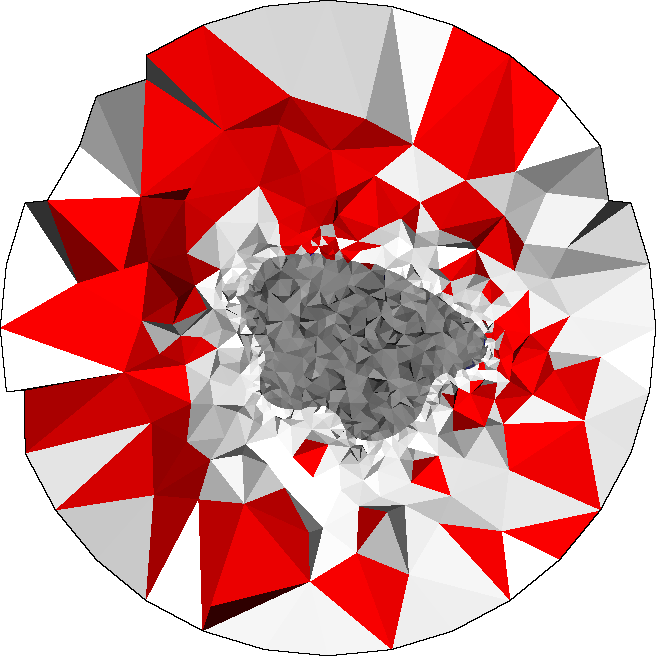} 
& \includegraphics[width=0.40\textwidth]{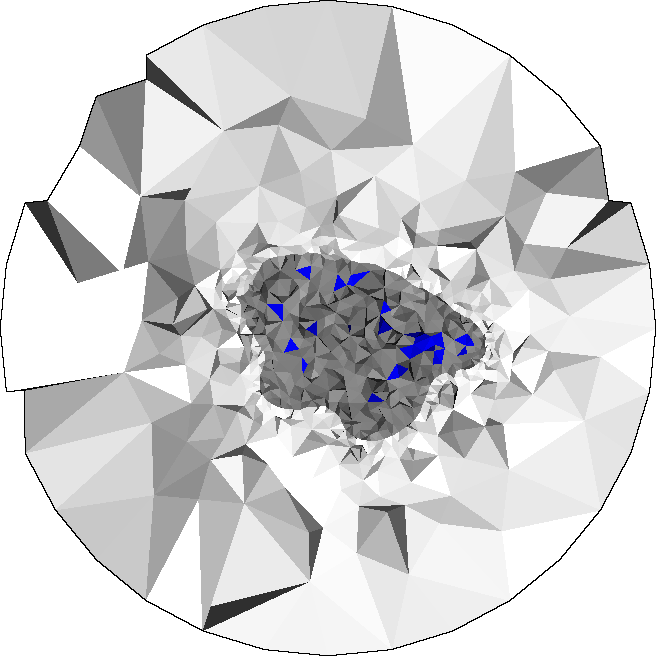} 
\end{tabular}
\end{center}
\caption{A cut-away of the $3$D mesh surrounding Fasciculin-$1$.  The colors indicate the
distribution of marked elements using the \emph{global} marking strategy 
with either the energy-based indicator (left) or the goal-oriented indicator (right).  Red and blue are marked elements in the
solvent and solute subdomains, respectively.}
\label{fig:1fas-global-mark}
\end{figure}

The reason for the poor performance can be explained by looking at Fig.~\ref{fig:1fas-global-mark}.
The images are cut-aways of the $3$D  mesh, colored to indicate the
distribution of marked elements. The left image shows elements selected by the
global marking strategy using the energy-based indicator. While the image
on the right uses the goal-oriented indicator.
The white elements are unmarked elements in the solvent subdomain and the gray
elements are unmarked elements in the solute subdomain.  Elements colored red
are marked solvent elements, while blue elements are marked solute elements.  
Notice for the energy-based indicator only elements in the solvent subdomain
are marked.  This indicates that the reaction
potential in the solute domain is relatively
well approximated compared to the solution in the solvent subdomain. The
distribution of the elements marked using a goal-oriented
indicator is focused on a few locations in the inner subdomain around
the solute atoms. 
This explains the poor convergence for the goal-oriented indicator 
in Fig.~\ref{fig:global-sfe-1fas}. Since the strategy does not indicate there
is error in the solvent domain, no refinement takes place.

\begin{figure}
\begin{center}
\includegraphics[width=0.40\textwidth]{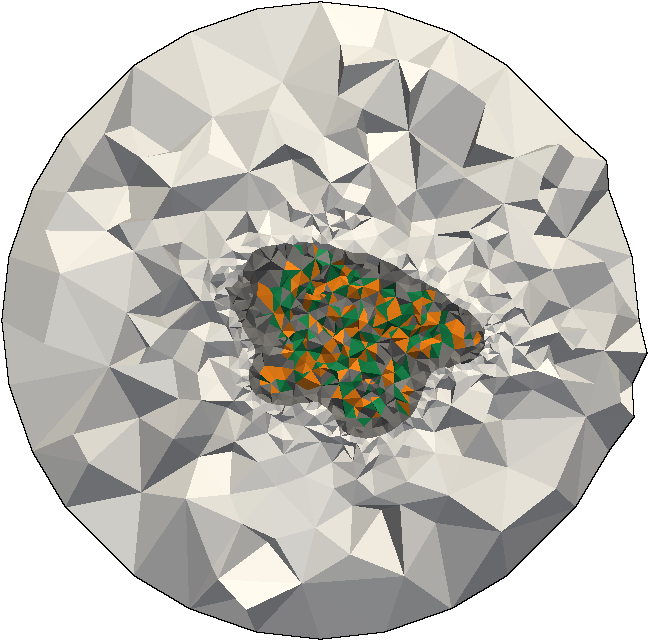}
\end{center}
\caption{A cut-away of the $3$D mesh surrounding Fasciculin-$1$.  The colors 
indicate the positive (green) and negative (orange) estimated elementwise 
contributions to the error in the solvation free energy (see eq.~\eqref{eqn:solv-err-eqn}).}
\label{fig:1fas-oscillations}
\end{figure}

To gain further insight into why the global goal-oriented strategy only marks
elements in the solute domain, we refer the reader to
Fig.~\ref{fig:1fas-oscillations}.
This image shows a cut-away with elements colored by their approximate
\emph{signed contribution} to the error in the goal. 
This is constructed using a quadratic approximation of the dual and the
signed elementwise contributions from Eq.~\eqref{eqn:solv-err-eqn}.
Note that because these are error contributions and not indicators, they take 
both positive and negative values. The element contributions range in value 
between $-0.6928$ and $1.2779$.  In the figure, positive element contributions 
greater than $0.005$ are colored green and negative contributions less 
than $-0.005$ are
colored orange. From the image it is clear that the contributions 
in the solute subdomain have relatively large magnitude, and they oscillate in sign.  The result of
this oscillation is that the majority of these contributions cancel when 
integrated over the entire solute domain.  However, the error
estimator in Eq.~\eqref{eqn:goal-quad-ind-first} uses the absolute value,
which results in an overestimation of the error attributed to the solute
domain. As a result, the solute domain is 
over refined, unless steps are taken to modify the marking (or error
estimation) strategy.

\subsubsection{Split Marking Strategy}
To improve on the convergence of goal-oriented refinement, a second marking strategy is employed.
This is a domain dependent marking strategy that attempts to spread the refinement over solvent and solute
regions of the domain.  The strategy relies on splitting the mesh into two subsets $\calT_l^s\subset\calT_l$ 
and $\calT_l^m\subset\calT_l$, where $\calT_l^s$ and $\calT_l^m$ contain the elements in the solvent and
solute domains respectively.  Stated concisely, the marking strategy is
\begin{equation}\label{eqn:pbe-split-mark}
\mbox{Mark all}\quad
\begin{array}{c}
K \in \calT_l^s \\
K \in \calT_l^m
\end{array}\quad
\mbox{such that}\quad
\begin{array}{c}
\eta_{K} > \gamma \; \underset{T \in \calT_l^s}{\max} \; \eta_T \\
\eta_{K} > \gamma \; \underset{T \in \calT_l^m}{\max} \; \eta_T
\end{array}
\end{equation}
where $\gamma \in (0,1)$.

\begin{figure}
\begin{tabular}{cc}
  {Energy-Based} & {Goal-Oriented} \\
  \includegraphics[width=0.40\textwidth]{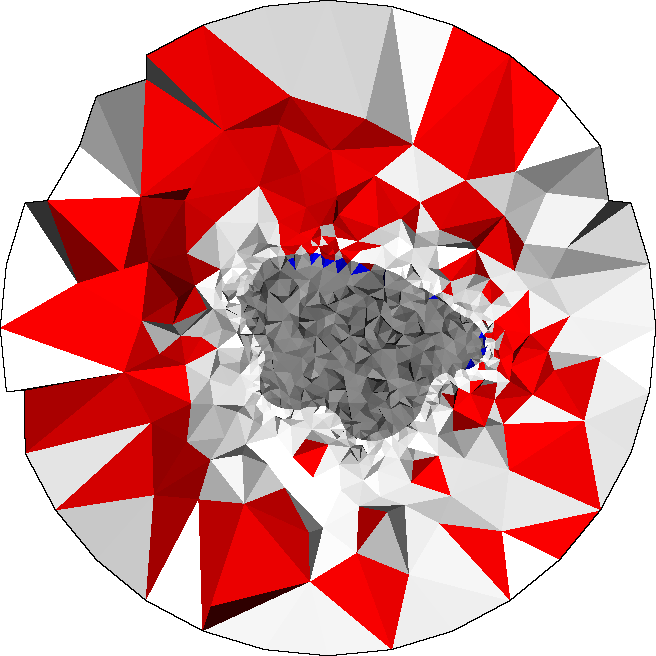}
& \includegraphics[width=0.40\textwidth]{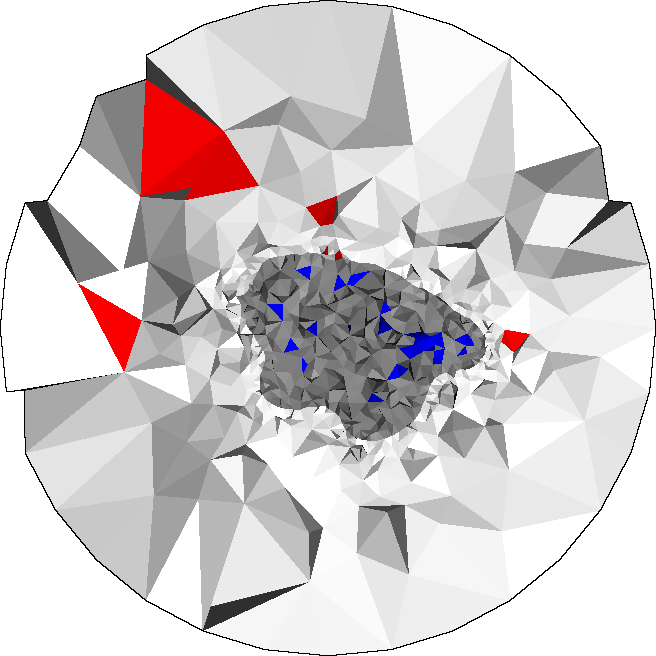}
\end{tabular}
\caption{A cut-away of the $3$D mesh surrounding Fasciculin-$1$.  The colors indicate the
distribution of marked elements using the \emph{split} marking strategy 
with either the energy-based indicator (left) or the goal-oriented indicator (right).  
Red and blue are marked elements in the
solvent and solute subdomains, respectively.}
\label{fig:1fas-split-mark}
\end{figure}

The split-marking strategy marks a significantly different group of elements, especially for the goal-oriented
indicators.  The color coding of Fig.~\ref{fig:1fas-split-mark} is the same as in Fig.~\ref{fig:1fas-global-mark}.  
Again the energy-based refinement selects elements primarily in the solvent subdomain. However, because split-marking
forces refinement in both subdomains, a few elements along the interface in the solute domain are also selected. In
contrast, split-marking using a goal-oriented indicator marks a few elements in the solvent subdomain, while also marking
elements surrounding the solute atoms. 

\begin{figure}
\begin{center}
\includegraphics[width=.45\textwidth]{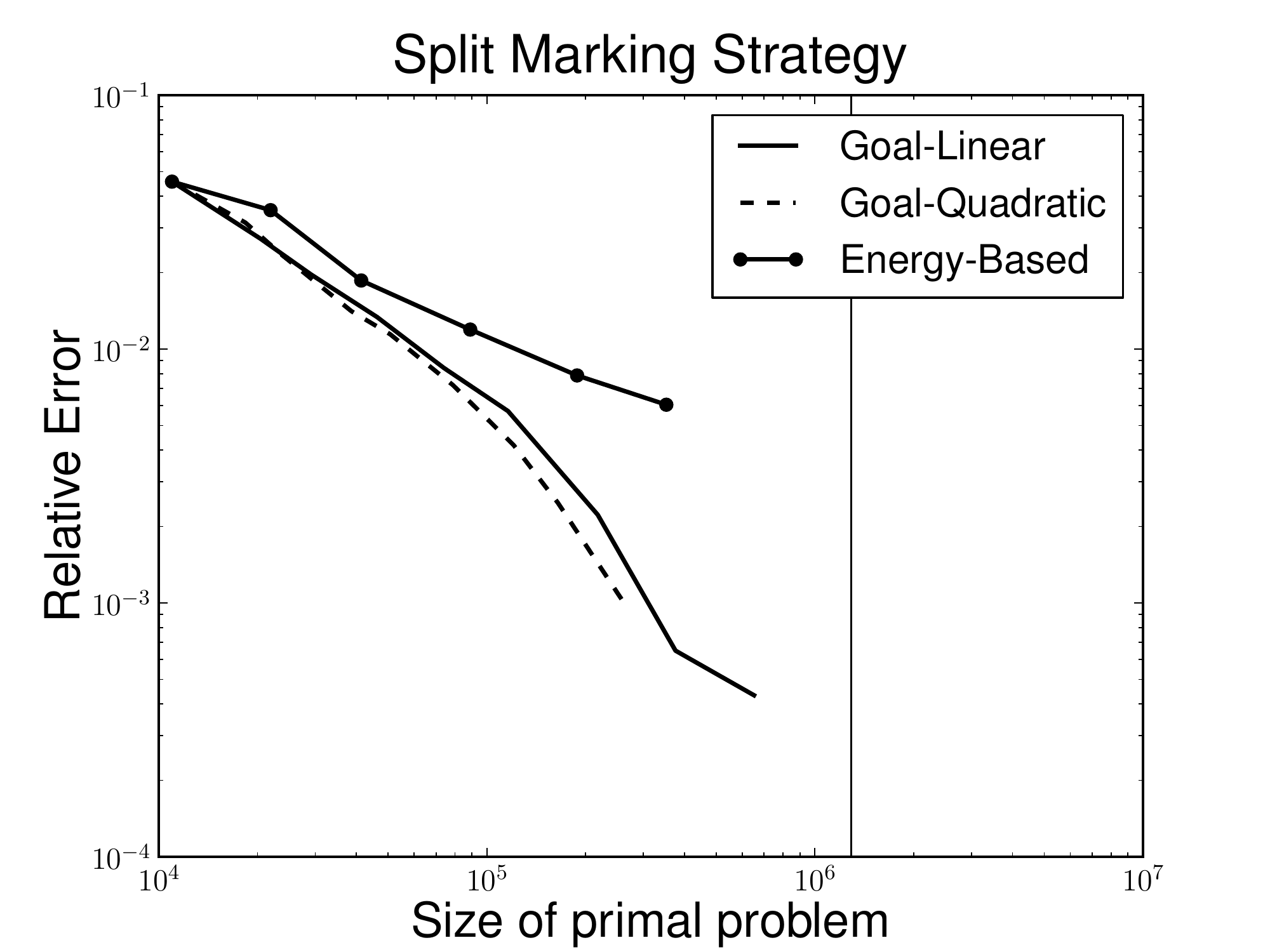}
\includegraphics[width=.45\textwidth]{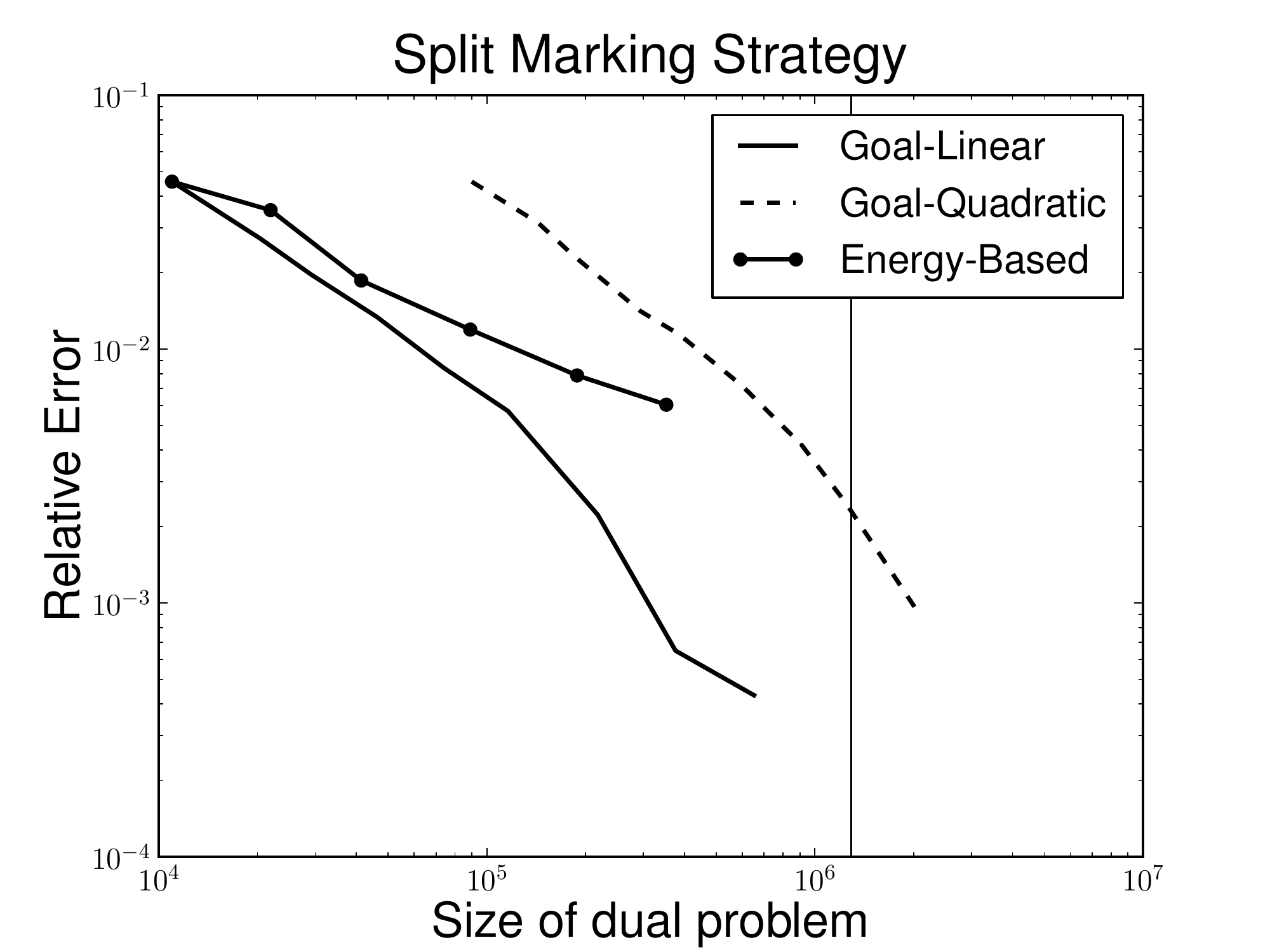}
\end{center}
\caption{Convergence of the solvation free energy for Fasciculin-$1$ with both goal-oriented and energy-based
indicators using the split marking strategy. The convergence measured against the size of the primal/dual problem
is on the left/right.}
\label{fig:split-sfe-1fas}
\end{figure}

Figure~\ref{fig:split-sfe-1fas} shows the relative error of the solvation free energy in Fasciculin-$1$ as
a function of problem size for several indicators using the split-marking strategy.  
In the figure, there are two plots. 
The plot on the left shows the relative error in the solvation free energy as a function of the
number of unknowns in the primal problem.  
The second plot, with the exception of the energy-based refinement strategy,
shows the relative error as a function of the size of the dual problem. For energy-based
refinement, since no dual is needed, the horizontal axis is the size of the primal problem.

In the figure, the split marking strategy using the energy-based indicator
converges steadily. This is similar to the convergence of energy-based
refinement using global marking.  On the other hand, there is a dramatic
improvement in the convergence of both the linear and quadratic goal-oriented refinement techniques,
with the rate even increasing slightly as the size of the problem increases.
Compare this with the poor results for goal-oriented refinement with the
global marking strategy from Fig.~\ref{fig:global-sfe-1fas}.
The improvement comes from the split marking strategy explicitly taking into
account the error in the solvent
and trying to control where it is large.

The second plot in Fig.~\ref{fig:split-sfe-1fas} shows that
solving the dual problem using piecewise quadratic elements has substantial
additional cost. Although, it is likely that this
cost would be mitigated by the additional work needed when solving the
nonlinear RPBE (see Eq.~\ref{eqn:rpbe}).  In contrast, the goal-oriented
strategy using an indicator constructed from a linear dual problem (and split marking) does not suffer from the
same problem, and is the most efficient method when the total cost is taken
into account.

\subsection{Performance of Preconditioners}
As was discussed in section~\ref{sec:prec}, classical multigrid (MG) 
preconditioners perform best in the uniform refinement setting, where there is
a rapid geometric growth in the number of unknowns as the mesh is refined.  In
the local refinement setting, this growth is much slower, and is frequently
subgeometric when the refinement is concentrated in the neighborhood of a
low dimensional feature (e.g., a point or a line).  As a result, the 
per-iteration complexity of classical MG may fail to scale linearly (or
suffer from a large scaling constant) as the number of unknowns increases.  In
contrast, many local MG preconditioners do not have the same restriction, and
maintain optimal per-iteration complexity in both the local and uniform 
refinement settings.

\begin{figure}
\begin{center}
\includegraphics[width=.75\textwidth]{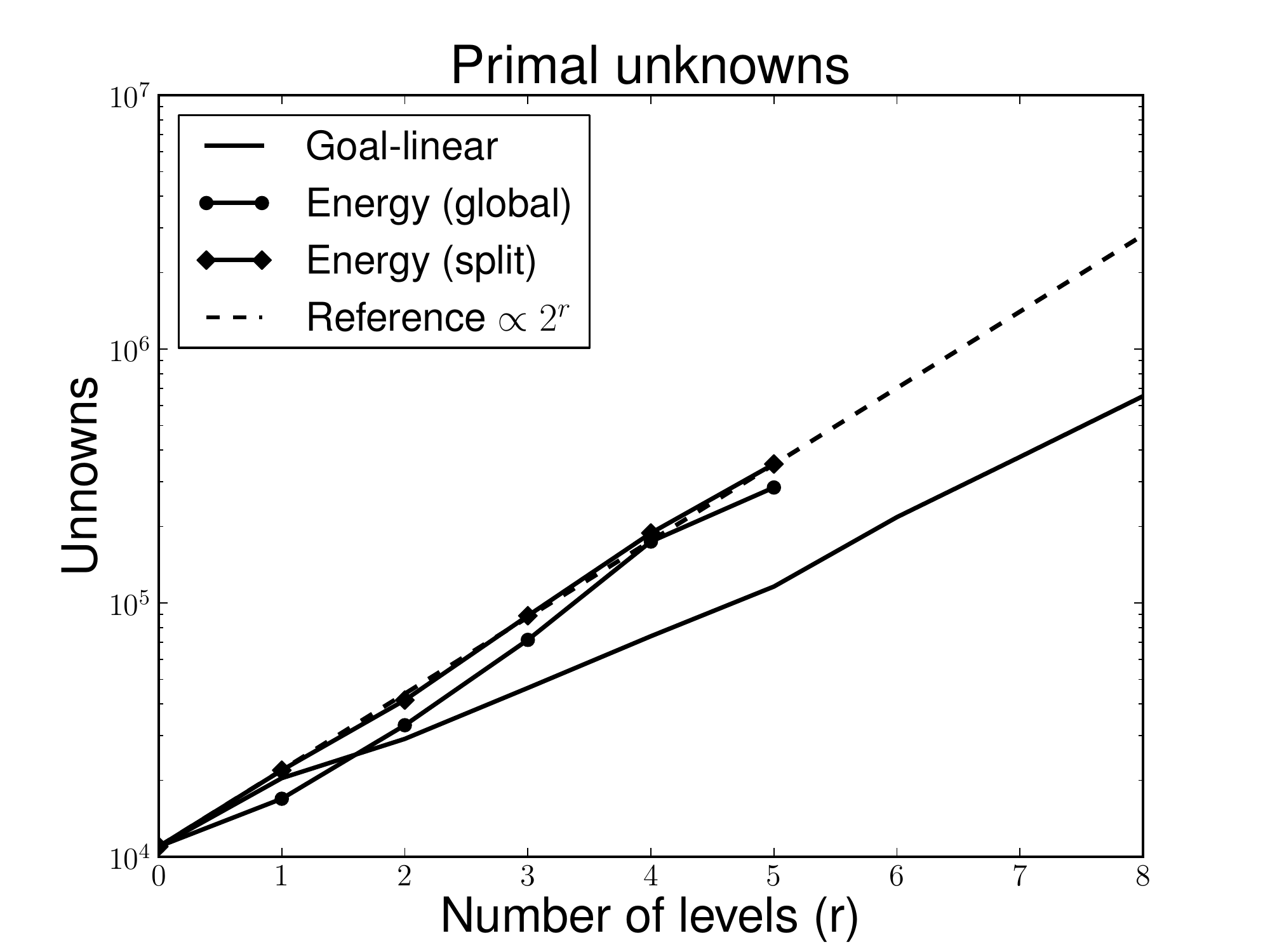}
\end{center}
\caption{Number of primal unknowns as a function of the number of levels 
of refinement. The reference line grows proportional to $2^r$.}
\label{fig:adapt-unk-growth}
\end{figure}
In Figure~\ref{fig:adapt-unk-growth}, the number of primal unknowns is shown as
a function of the refinement level for three different refinement strategies.
For the energy-based marking/refinement strategy, the growth in the number
of unknowns is geometric, but subuniform.  For goal-oriented refinement  
the growth is even slower, but still geometric.  The dashed
reference line shows a geometric growth rate proportional to $2^r$.

\begin{figure}
\begin{center}
\includegraphics[width=.75\textwidth]{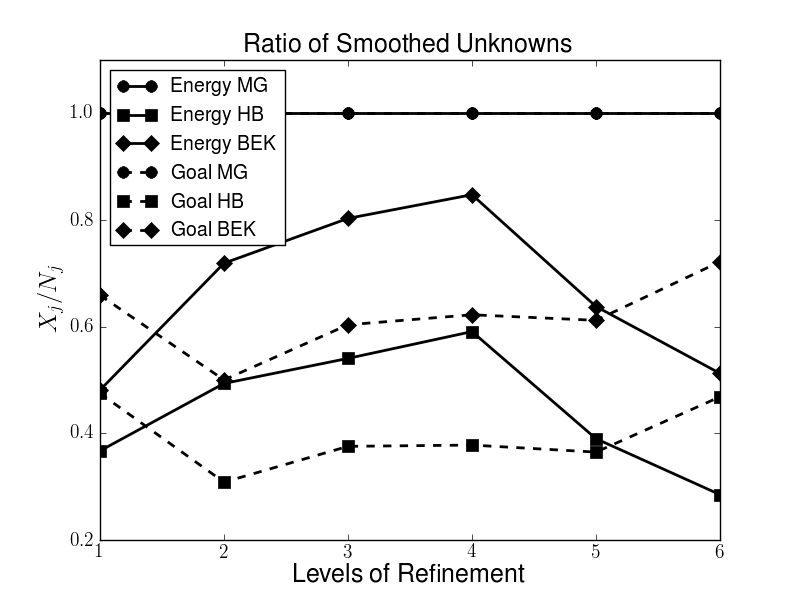}
\end{center}
\caption{Ratio of the number smoothed unknowns, $X_j$, to total
unknowns, $N_j$, on each level for different multilevel preconditioners}
\label{fig:onefas-ratio-dof}
\end{figure}
As a measure of the relative locality of each preconditioner, we compute
the ratio of the number of unknowns processed by the smoother, shown in
Figure~\ref{fig:onefas-ratio-dof}.  Classical MG is a global method, so
the ratio of smoothed to total unknowns is 1 in both the energy based and
goal oriented refinement cases.  As expected, the BEK preconditioner
consistently smooths more unknowns than HB, but fewer than MG.  This is
because the set processed by the BEK is a superset of the set processed
by HB.

\begin{figure}
\begin{center}
\includegraphics[width=.75\textwidth]{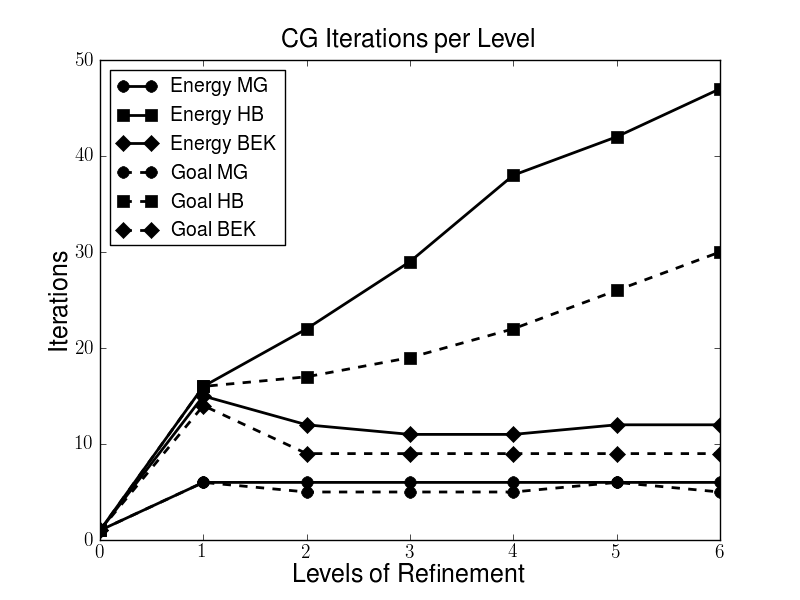}
\end{center}
\caption{Conjugate gradient method iteration counts for the preconditioners used.}
\label{fig:onefas-iters-dof}
\end{figure}
By design, the local multilevel preconditioners considered here maintain an 
optimal per-iteration complexity in both the local and global refinement 
settings. The primary challenge for these preconditioners is achieving a
bounded condition number, independent of problem size. It can be shown that 
for a given tolerance, the number of conjugate gradient (CG) iterations can be
bounded by a function of the condition number.  Hence, if the condition
number is bounded, so is the number of CG iterations. In 
Figure~\ref{fig:onefas-iters-dof}, we report the number of CG iterations as
a function of refinement level for MG, BEK, and HB.  As 
predicted by theory, application of the HB preconditioner leads to a 
slow growth in the number of CG iterations as the mesh is refined regardless
of indicator type. In contrast, for both goal-oriented and energy-based refinement
the iteration count is bounded for both the classical MG and 
BEK preconditioners.  The iteration counts for BEK are modestly
higher than classical MG, but the work per level is reduced, since BEK
smooths only a fraction of the unknowns smoothed by classical MG.  For
this reason, BEK is a compelling alternative to classical MG and HB.

\section{Conclusion}
\label{sec:conc}
In this article, we developed goal-oriented error indicators for accurate 
computation of the solvation free energy from solutions of the 
regularized Poisson-Boltzmann equation.  We found that due to oscillations
and imbalanced cancellation in the error contributions, global marking
strategies based on goal-oriented 
error indicators were not viable for driving adaptive mesh refinement. 
To address this problem, we developed a split marking strategy based on 
considering each subdomain individually.  In numerical experiments, we
calculated the solvation free energy for
a 921-atom Fasciculin-1 protein. Through these experiments, we showed
that the new marking strategy, combined with goal-oriented refinement, is more 
efficient than energy-based refinement in the context of solvation free energy
calculations.

The use of adaptive mesh refinement puts a greater burden on the 
preconditioner to maintain optimal runtime efficiency. To address this issue, 
we investigated the use of local multigrid methods, which have a lower 
per-iteration complexity compared to classical global multigrid. In particular,
the BEK variant proved to be a compelling alternative to classical 
multigrid since it has optimal per-iteration complexity, while still 
maintaining a bounded iteration count as the mesh is refined.
The result is an iterative solver with an optimal overall complexity,
scaling linearly with problem size.

\section{Acknowledgments}
\label{sec:ack}
BA was supported in part by NSF Award~1016190.  SB was supported in
part by NSF Award~0830578. EC was supported in part by a University of
Illinois CSE Fellowship and by a DOE Office of Science ASCR-UQ effort
at Sandia National Laboratory under contract DE-AC04-94AL85000.
MH was supported in part by NSF Awards~0715146 and 0915220, and
by DOD/DTRA Award HDTRA-09-1-0036.
Sandia National Laboratories is a multi-program laboratory managed and 
operated by Sandia Corporation, a wholly owned subsidiary of Lockheed Martin 
Corporation, for the U.S. Department of Energy's National Nuclear Security 
Administration under contract DE-AC04-94AL85000.

\bibliographystyle{abbrv}
\bibliography{refs}

\vspace*{0.5cm}

\end{document}